\documentclass[11pt,a4paper]{article}

\usepackage{amsmath}
\usepackage{amssymb}
\usepackage{graphicx}
\newtheorem{lem}{Lemma.}[section]

\newtheorem{prop}{Proposition.}[section]
\newtheorem{asu}{Assumption.}[section]

\author{Atsushi Iwasaki\thanks{Fukuoka Institute of Technology\ \  e-mail: a-iwasaki@fit.ac.jp}}
\title{Deriving the Variance of the Discrete Fourier Transform Test\\
Using Parseval's Theorem}
\begin{document}

\maketitle
\abstract{
The discrete Fourier transform test is a randomness test included in NIST SP800-22.
However, the variance of the test statistic is smaller than expected and the theoretical value of the variance is not known.
Hitherto, the mechanism explaining why the former variance is smaller than expected has been qualitatively explained based on Parseval's theorem.
In this paper, we explore this quantitatively and derive the variance using Parseval's theorem under particular assumptions. Numerical experiments are then used to show that this derived variance is robust. 
}

\section{Introduction}
Random number sequences are used in many fields.
These sequences play a particularly important role in information security including cryptography because high ``randomness'' is required in such fields.
Thus, evaluation of random number sequences and their generators is indispensable.

Randomness tests are one of the most fundamental evaluation methods in this respect.
These are hypothesis tests, and the null is that the given sequence is truly random.
Randomness tests do not require information about the generator which resulted in the given sequence.
Thus, such tests can be widely used without regard to the generators.
There are many kinds of randomness tests and some test sets have been proposed.
NIST SP800-22 \cite{NIST} is one of the most well-known test sets; the first version was published in 2001 and revision 1a, published in 2010, is currently the most recent variant.
Revision 1a consists of 188 test items which can be grouped according to 15 types.

A discrete Fourier transform test (DFTT) is one of the randomness tests included in NIST SP800-22.
The algorithm for the DFTT included in the first version is as follows:
\begin{enumerate}
\item Input an $n$-bit sequence $X$. Here, each bit of $X$ is 0 or 1.
\item Convert each bit $x$ to $2x-1$ (each bit becomes 1 or -1).
\item Perform DFT to obtain the Fourier spectrum series 
\begin{align}
|f_0|,|f_1|,\cdots,|f_{\left\lfloor\frac{n}{2}\right\rfloor-1}|,
\end{align}
where
\begin{align}
|f_j|:=\sqrt{\left(\sum_{k=0}^{n-1}x_k\cos\left(\frac{2\pi kj}{n}\right)\right)^2+\left(\sum_{k=0}^{n-1}x_k\sin\left(\frac{2\pi kj}{n}\right)\right)^2},
\end{align}
with $x_k$ being the $(k+1)$-th bit of $X$.
\item Count the elements of $\{|f_0|,|f_1|,\cdots,|f_{\left\lfloor\frac{n}{2}\right\rfloor-1}\}$ satisfying $|f_i|<\sqrt{3n}$.
Let $N_1$ be such a number.
\item Compute $d$ defined by
\begin{align}
d=\frac{N_1-0.95\frac{n}{2}}{\sqrt{(0.95)(0.05)\frac{n}{2}}}.
\end{align}
\item Compute p-value $p$ defined by
\begin{align}
p=\text{erfc}\left(\frac{|d|}{\sqrt{2}}\right),
\end{align}
where erfc denotes the complementary error function.
\end{enumerate}
For $j=1,2,\cdots,n-1$, $\frac{2}{n}|f_j|^2$ follows the $\chi^2$-distribution with 2 degrees of freedom as $n\to\infty$.
If $Z$ is a stochastic variable and $\frac{2}{n}|Z|^2$ follows a $\chi^2$-distribution with 2 degrees of freedom, then
\begin{align}
\text{Prob}\{\ |Z|<\sqrt{3n}\ \}\sim 0.95.
\end{align}
Thus, it is expected that $N_1$ approximately follows a normal distribution with an average and variance of $0.95\frac{n}{2}$ and $(0.95)(0.05)\frac{n}{2}$, respectively.
If this holds, $d$ and $p$ approximately follow a standard normal distribution and a uniform distribution over the interval $[0,1]$, respectively.

However, the following problems were revealed by \cite{Kim-Umeno-Hasegawa,Hamano}.
\begin{itemize}
\item The threshold $\sqrt{3n}$ is an approximation and it is not accurate enough for practical use.
\item The variance of $N_1$ observed numerically is far smaller than expected.
\end{itemize}
To address the first problem, Kim et al. proposed that the threshold value should be changed from $\sqrt{3n}$ to $\sqrt{-n\log(0.05)}$ \cite{Kim-Umeno-Hasegawa}.
If $Z$ is a stochastic variable and $\frac{2}{n}|Z|^2$ follows a $\chi^2$-distribution with 2 degrees of freedom,
\begin{align}
\text{Prob}\{\ |Z|<\sqrt{-n\log(0.05)}\ \}=0.95\ .
\end{align}
Then, we can state that the proposed threshold is appropriate. The DFTT included in revision 1a uses this proposed value.

To address the second problem, Kim et al. also proposed to use $(0.95)(0.05)\frac{n}{4}$ as an estimate of the variance of $N_1$ \cite{Kim-Umeno-Hasegawa}.
In other words, because the variable $d$ should follow a standard normal distribution, they proposed to change the computation of $d$ to
\begin{align}
d=\frac{N_1-0.95\frac{n}{2}}{\sqrt{(0.95)(0.05)\frac{n}{4}}}.
\end{align}
Subsequently, Hamano et al. proposed $(0.95)(0.05)\frac{n}{2}\times0.528\simeq(0.95)(0.05)\frac{n}{3.7879}$ as the variance of $N_1$ \cite{Hamano-Yamamoto} and Pareschi et al. proposed $(0.95)(0.05)\frac{n}{3.8}$ as that \cite{Pareschi}.
These proposed values were experimentally derived but the theoretical value of this variance has not been derived.

Yamamoto et al. qualitatively explained why the variance of $N_1$ is smaller than\\ $(0.95)(0.05)\frac{n}{2}$ based on Parseval's theorem \cite{Yamamoto-Kaneko}.
We believe that the explanation is persuasive.
In this paper, we derive the variance of $N_1$ based on Parseval's theorem with some assumptions.
Although some alternative tests to DFTT have been proposed \cite{Okada,Iwasaki-SCIS}, these are not dealt with herein, because the DFTT is now preferred and revising a parameter is easier than revising the algorithms.

The remainder of the paper is organized as follows.
In section 2, we introduce the qualitative explanation based on Parseval's theorem which has hitherto been documented in the literature.
In section 3, we theoretically derive the value of the variance with Parseval's theorem and some assumptions.
In section 4, we perform some experiments and discuss the validity of the result in Section 3.
Finally, in section 5 we offer conclusions.

\section{Qualitative explanation based on Parseval's theorem}
In this section, we briefly introduce the qualitative explanation by Yamamoto et al. This also serves as preparation for the next section.

By Parseval's theorem, it follows that
\begin{align}
\sum_{j=0}^{n-1}|f_j|^2=n\sum_{k=0}^{n-1}|x_k|^2.
\end{align}
For $k=0,1,\cdots,n-1$, $|x_k|^2=1$ because $x_k=\pm1$.
Then,
\begin{align}
\sum_{j=0}^{n-1}|f_j|^2=n^2.
\end{align}
By the symmetry of the Fourier transformation, 
\begin{align}
|f_{n-j}|=|f_j|
\end{align}
for arbitrary $j$.
Then,
\begin{align}
\sum_{j=0}^{\left\lfloor\frac{n}{2}\right\rfloor-1}|f_j|^2=
\begin{cases}
\frac{n^2}{2}+\frac{|f_0|^2}{2}-\frac{|f_{\frac{n}{2}}|^2}{2}\ \ \ \ &(n\text{:even})\\
\frac{n^2}{2}+\frac{|f_0|^2}{2}\ \ &(n\text{:odd})
\end{cases}
.\label{restriction}
\end{align}
From \cite{Hirose}, for $j=1,2,\cdots,n-1$, we have
\begin{align}
V\left[|f_j|^2\right]=n^2-2n, \label{hirose1}
\end{align}
where $V[Z]$ is the variance of variable $Z$.
For $j=0$, it follows that
\begin{align}
V\left[|f_0|^2\right]&=V\left[\left(\sum_{k=0}^{n-1}x_k\right)^2\right]\\
&=\mathbb{E}\left[\left(\sum_{k=0}^{n-1}x_k\right)^4\right]-\mathbb{E}\left[\left(\sum_{k=0}^{n-1}x_k\right)^2\right]^2\\
&=\sum_{k_1=0}^{n-1}\sum_{k_2=0}^{n-1}\sum_{k_3=0}^{n-1}\sum_{k_4=0}^{n-1}\mathbb{E}\left[x_{k_1}x_{k_2}x_{k_3}x_{k_4}\right]-\left\{\sum_{k_1=0}^{n-1}\sum_{k_2=0}^{n-1}\mathbb{E}\left[x_{k_1}x_{k_2}\right]\right\}^2\\
&=2n^2-2n,\label{hirose2}
\end{align}
where $\mathbb{E}[\cdot]$ is $\frac{1}{2^n}\sum_{X\in\{-1,1\}^n}(\cdot)$\ .
Using (\ref{hirose1}) and (\ref{hirose2}),
\begin{align}
V\left[\sum_{j=0}^{\left\lfloor\frac{n}{2}\right\rfloor-1}|f_j|^2\right]\leq& \max\left\{4V\left[\frac{|f_0|^2}{2}\right],4V\left[\frac{|f_{\left\lfloor\frac{n}{2}\right\rfloor}|^2}{2}\right]\right\}\\
=&2n^2-2n.\label{varianceofenergy}
\end{align} 
If $|f_0|, |f_1|, \cdots, |f_{\left\lfloor\frac{n}{2}\right\rfloor-1}|$ are mutually independent, then
\begin{align}
V\left[\sum_{j=0}^{\left\lfloor\frac{n}{2}\right\rfloor-1}|f_j|^2\right]
=&\sum_{j=0}^{\left\lfloor\frac{n}{2}\right\rfloor-1}V\Big[|f_j|^2\Big]\\
=&\left\lfloor\frac{n}{2}\right\rfloor(n^2-2n)+n^2.\label{varianceofenergywithassumption}
\end{align} 
Comparing (\ref{varianceofenergy}) and (\ref{varianceofenergywithassumption}), we conclude that the energy $\sum_{j=0}^{\left\lfloor\frac{n}{2}\right\rfloor-1}|f_j|^2$ is restricted to a narrow area in $\mathbb{R}$ by the dependency among $|f_0|, |f_1|, \cdots, |f_{\left\lfloor\frac{n}{2}\right\rfloor-1}|$.

Yamamoto et al. claimed that this restriction explains why $V[N_1]$ is smaller than $(0.95)(0.05)\frac{n}{2}$.
If some elements of $\{|f_0|, |f_1|, \cdots, |f_{\left\lfloor\frac{n}{2}\right\rfloor-1}|\}$ take large values, then some other elements must take small values to maintain the restriction.
As a result, the probability that $N_1$ takes a value far from the average is forced to be small, i.e., $V[N_1]$ is smaller than $(0.95)(0.05)\frac{n}{2}$.

\section{Quantitative analysis based on Parseval's theorem}
From our perspective, the qualitative explanation introduced in the previous section is persuasive.
The restriction by Parseval's theorem will be critical to deriving $V[N_1]$.
Thus, in this section, we theoretically derive $V[N_1]$ with Parseval's theorem under particular assumptions.
To simplify discussion we assume that $n$ is even, but this is not essential.
We define $m$ as $\frac{n}{2}$.

\subsection{Assumptions}
If we know the analytical form of the joint probability density function depending on \\$\left(|f_0|, |f_2|, \cdots, |f_{m-1}|\right)$, we can derive $V[N_1]$.
However, the form is not known.
Thus, we need to assume a certain analytical form and in this subsection we discuss which such form is appropriate. 
First, we introduce the following proposition.
\begin{prop}
\label{prop-3}
Let $R$ be an arbitrary natural number.
Then, the joint probability density function depending on $\left(\frac{1}{m}|f_1|^2, \frac{1}{m}|f_2|^2, \cdots, \frac{1}{m}|f_R|^2\right)$ converges to
\begin{align}
\frac{1}{2^R}\exp\left(-\frac{\sum_{j=1}^R|f_j|^2}{2m}\right)
\label{pdf}
\end{align}
as $m\to \infty$.
\end{prop}
In Proposition \ref{prop-3}, $R$ is a fixed value.
As an analogy, let us assume that the joint probability density function depending on $\left(\frac{1}{m}|f_0|^2, \frac{1}{m}|f_2|^2, \cdots, \frac{1}{m}|f_{m-1}|^2\right)$ converges to 
\begin{align}
\frac{1}{2^m}\exp\left(-\frac{\sum_{j=0}^{m-1}|f_j|^2}{2m}\right).
\label{pdf2}
\end{align} 
However, there are two problems:
\begin{itemize}
\item Even as $n\to\infty$, i.e., $m\to\infty$, $\frac{1}{m}|f_0|^2$ does not follow a $\chi^2$-distribution with 2 degrees of freedom, which is the limit distribution of $\frac{1}{m}|f_j|^2,\ (j\ne 0)$.
\item The energy is restricted to a narrow area by (\ref{restriction}).
\end{itemize}
To address these problems, we assume the following:
\begin{itemize}
\item For sufficiently large $n$, the value of $V[N_1]$ will be virtually preserved even if $\frac{1}{m}|f_0|^2$ are replaced with a variable following a $\chi^2$-distribution with 2 degrees of freedom.
Thus, we assume that $\frac{1}{m}|f_0|^2$ follows a $\chi^2$-distribution with 2 degrees of freedom.
\item The energy $\sum_{j=0}^{m-1}|f_j|^2$ is restricted by (\ref{restriction}), and it is not constant.
However, by (\ref{varianceofenergy}), we can state that the standard deviation of the energy is sufficiently smaller than the average of that energy. 
Thus, we assume that
\begin{align}
\sum_{j=0}^{m-1}|f_j|^2=2m^2.\label{constantrestriction}
\end{align}
\end{itemize}
Based on the above two assumptions and (\ref{pdf2}), we further assume the following.
\begin{asu}
\label{assumption}
In $\mathbb{R}^m$, $\left(|f_0|^2, |f_2|^2, \cdots, |f_{m-1}|^2\right)$ is restricted on the surface defined by (\ref{constantrestriction}), and it is uniformly distributed in the area where $|f_j|^2\geq0\ j=0,1,2,\cdots,m-1$ on the surface.
\end{asu}

\subsection{Analysis using polar coordinates}

To derive $V\left[N_1\right]$, we need probability density functions depending on $|f_i|^2$ for $i=0,1,\cdots,m-1$ and joint probability density functions depending on $\left(|f_i|^2,|f_j|^2\right)$ for $i\ne j$.
In this subsection, we derive these different density functions under Assumption \ref{assumption}.

Under Assumption \ref{assumption}, $\Bar{p}$ which is the joint probability density function depending on $\left(|f_0|^2, |f_1|^2, \cdots, |f_{m-1}|^2\right)$ is as follows:
\begin{align}
\Bar{p}(|f_0|^2, |f_1|^2, \cdots, |f_{m-1}|^2)=
\begin{cases}
C_1&(|f_j|^2\geq0\ (j=0,1,\cdots,m-1)\text{ and }\sum_{j=0}^{m-1}|f_j|^2=2m^2)\\ 
0&(\text{otherwise})
\end{cases}
.\label{uniformdistribution}
\end{align}
We use polar coordinates,
\begin{align}
 |f_0| &= r \cos\theta_1, \label{X_0}\\
      |f_1| &= r \sin\theta_1 \cos\theta_2, \\
      |f_2|  &= r\sin\theta_1 \sin\theta_2 \cos\theta_3, \\
              &\vdots \\
      |f_{m-2}| &= r \sin\theta_1 \cdots \sin\theta_{m-2} \cos \theta_{m-1}, \\
      |f_{m-1}| &= r \sin\theta_1 \cdots \sin\theta_{m-2} \sin \theta_{m-1},
\end{align}
where $r\geq 0$, $0\leq \theta_i\leq \pi$ ($i=1,2,\cdots,m-2$), and $0\leq \theta_{m-1}< 2\pi$.
The Jacobian for the change in coordinates is
\begin{align}
r^{m-1} \prod_{i=1}^{m-1}\left(\sin\theta_i\right)^{m-i-1}.
\end{align}
For $j=0,1,\cdots,m-1$, it follows that
\begin{align}
d|f_j|^2=2|f_j|d|f_j|.
\end{align}
Then, we have
\begin{align}
\prod_{j=0}^{m-1}d|f_j|^2
=&2^m\prod_{j=0}^{m-1}|f_j|d|f_j|.\\
=&2^mr^{m-1}\left(\prod_{j=0}^{m-1}|f_j|\right) \left(\prod_{i=1}^{m-1}\left(\sin\theta_i\right)^{m-i-1}\right) drd\theta_1d\theta_2\cdots d\theta_{m-1}\\
=&2^mr^{2m-1}\left(\prod_{i=1}^{m-1}\left(\sin\theta_i\right)^{2m-2i-1}\cos \theta_i\right) drd\theta_1d\theta_2\cdots d\theta_{m-1}.
\label{henkan}
\end{align}

Let $p_j$ be a probability density function depending on $|f_j|^2$ and $p_{(i,j)}$ be a joint probability density function depending on $(|f_i|^2,|f_j|^2)$.
First, we derive $p_0$.
By (\ref{X_0}), $|f_0|^2=r^2\cos^2\theta_1$.
Then, by (\ref{uniformdistribution}) and (\ref{henkan}),
\begin{align}
p_0(|f_0|^2)d|f_0|^2=
\begin{cases}
C_2 \left(\sin\theta_1\right)^{2m-3}\cos \theta_1 d\theta_1\ \ &(0\leq |f_0|^2\leq 2m^2)\\
0 &(\text{otherwise})
\end{cases}
,
\end{align}
and
\begin{align}
d|f_0|^2=C_3\sin \theta_1 \cos \theta_1d\theta_1,
\end{align}
where $C_2$ and $C_3$ are constants. 
Then,
\begin{align}
p_0(|f_0|^2)=
\begin{cases}
\frac{m-1}{2m^2} \left(1-\frac{|f_0|^2}{2m^2}\right)^{m-2}\ \ &(0\leq |f_0|^2\leq 2m^2)\\
0 &(\text{otherwise})
\end{cases}
.
\end{align}
Because $|f_0|^2,|f_1|^2,\cdots|f_{m-1}|^2$ follow the same distribution under Assumption \ref{assumption}, for $j=0,1,\cdots,m-1$,
\begin{align}
p_j(|f_j|^2)=
\begin{cases}
\frac{m-1}{2m^2} \left(1-\frac{|f_j|^2}{2m^2}\right)^{m-2}\ \ &(0\leq |f_j|^2\leq 2m^2)\\
0 &(\text{otherwise})
\end{cases}
.\label{syuuhen}
\end{align}

By (\ref{uniformdistribution}), $p_{(i,j)}$ depends only on $|f_i|^2+|f_j|^2$.
Then, by (\ref{syuuhen}),
\begin{align}
&p_{(i,j)}(|f_i|^2,|f_j|^2)\nonumber\\
=&
\begin{cases}
\frac{(m-1)(m-2)}{(2m^2)^2} \left(1-\frac{|f_i|^2+|f_j|^2}{2m^2}\right)^{m-3}&(|f_i|^2\geq0, |f_j|^2\geq0\text{ and }|f_i|^2+|f_j|^2\leq 2m^2)\\
0 &(\text{otherwise})
\end{cases}
.\label{niko}
\end{align}

\subsection{Theoretical derivation of $V[N_1]$}
\begin{figure}[h]
\begin{center}
\includegraphics[width=10cm]{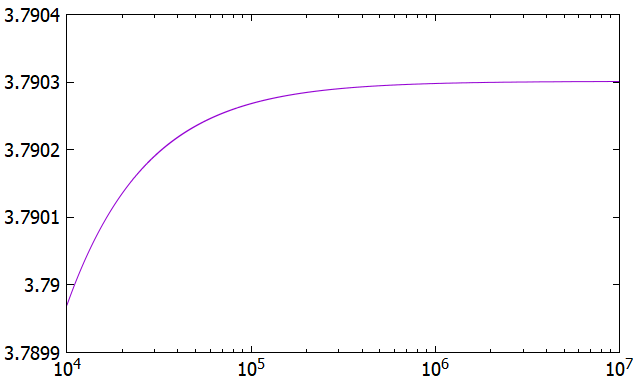}
\end{center}
\begin{picture}(0,0)
\put(240,10){\Large $m$}
\put(60,110){\rotatebox{90}{\Large $a$}}
\end{picture}
\caption{Plot of (\ref{riron})}
\label{rironkyokusen}
\end{figure}

In this subsection, we derive $V[N_1]$ using the results from the previous subsection.
For notational convenience, we introduce $T=\sqrt{-n\log(0.05)}$.
With this $T$, we define
\begin{align}
F_j:=
\begin{cases}
1\ \ &(|f_j|\leq T)\\
0   &(|f_j|>T)
\end{cases}
,
\end{align}
for $j=0,1,\cdots,m-1$.
Then, the variance is decomposed as  
\begin{align}
V[N_1]=&V\left[\sum_{j=0}^{m-1}F_j\right]\\
=&\sum_{j=0}^{m-1}V\Big[F_j\Big]+\sum_{(i,j)|i\ne j}C\Big[F_i, F_j\Big]\sqrt{V\Big[F_i\Big]V\Big[F_j\Big]},
\label{tenkai}
\end{align}
where $C[X,Y]$ is the correlation coefficient between variables $X$ and $Y$.
By (\ref{syuuhen}) and (\ref{niko}), we have $V[F_i]$ and $C[F_i,F_j]$ as follows.
\begin{align}
&V\left[F_i\right]=\left(1-\frac{T^2}{2m^2}\right)^{m-1}-\left(1-\frac{T^2}{2m^2}\right)^{2m-2},\label{bunsan}\\
&C\left[F_i,F_j\right]=\frac{\left(1-\frac{T^2}{m^2}\right)^{m-1}-\left(1-\frac{T^2}{2m^2}\right)^{2m-2}}{\sqrt{V\left[F_i\right]V\left[F_j\right]}}.\label{soukan}
\end{align}
Substituting (\ref{bunsan}) and (\ref{soukan}) into (\ref{tenkai}), we arrive at
\begin{align}
V\left[N_1\right]
=m\times&\left\{\left(1-\frac{T^2}{2m^2}\right)^{m-1}-\left(1-\frac{T^2}{2m^2}\right)^{2m-2}\right\}\nonumber\\
&+m(m-1)\times\left\{\left(1-\frac{T^2}{m^2}\right)^{m-1}-\left(1-\frac{T^2}{2m^2}\right)^{2m-2}\right\}.
\label{rironbunsan}
\end{align}
Assuming that $V\left[N_1\right]=(0.95)(0.05)\frac{n}{a}$, we have
\begin{align}
a
=\frac{2\times0.05\times0.95}{\left(1-\frac{T^2}{2m^2}\right)^{m-1}-\left(1-\frac{T^2}{2m^2}\right)^{2m-2}+(m-1)\times\left\{\left(1-\frac{T^2}{m^2}\right)^{m-1}-\left(1-\frac{T^2}{2m^2}\right)^{2m-2}\right\}}.
\label{riron}
\end{align}
Figure \ref{rironkyokusen} shows that $a$ approximately converges to 3.7903 as $m\to\infty$.
Broadly speaking, this is consistent with previous studies.

\section{Numerical experiments}
In the previous section, we used Assumption \ref{assumption}.
It is obvious that the true distribution is never the distribution defined by Assumption \ref{assumption}.
 However, they are very close in some sense and so we expect that (\ref{rironbunsan}) and (\ref{riron}), which are derived based on Assumption \ref{assumption}, are appropriate.
In this section, we discuss the validity of (\ref{rironbunsan}) and (\ref{riron}).

\subsection{Experiment 1}
By (\ref{tenkai}) which holds for both distributions, if $a$ converges to a certain value as $m\to\infty$, then the limit value depends only on the leading terms of $V[F_i]$ and $C[F_i,F_j]$.
By (\ref{bunsan}), for Assumption \ref{assumption},
\begin{align}
\lim_{m\to\infty}V[F_j]=(0.05)(0.95).
\label{bunsanlim}
\end{align}
With the true distribution, $\frac{1}{m}|f_j|^2$ converges to the $\chi^2$-distribution with 2 degrees of freedom for all $j$ except $j=0$, and (\ref{bunsanlim}) also holds.
We do not need to consider the case that $j=0$, because $V[F_0]$ becomes relatively smaller than $\sum_{j=0}^{m-1}V\left[F_j\right]$ as  $m\to\infty$.

\begin{figure}[h]
\begin{center}
\includegraphics[width=10cm]{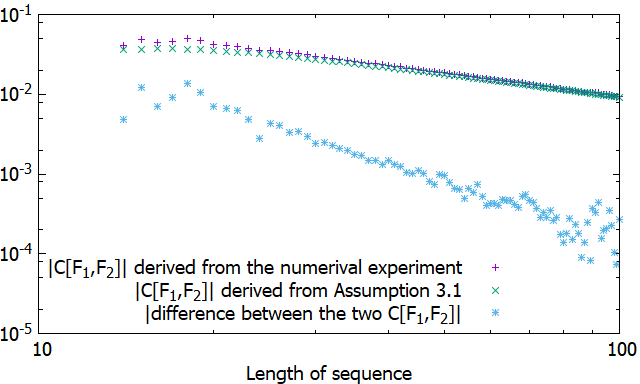}
\end{center}
\caption{Comparing $C[F_1,F_2]$}
\label{cor}
\end{figure}

Then, we should investigate the leading term of $C[F_i,F_j]$.
By the same argument, cases where $i=0$ and $j=0$ can be excepted.
We generated $10^8$ sequences using the Mersenne Twister \cite{MT}, performed DFT, and computed $C[F_1,F_2]$.
Figure \ref{cor} compares the $C[F_1,F_2]$ with that derived from Assumption \ref{assumption}.
The results corroborate that the leading term of $C[F_i,F_j]$ in both distributions is the same.

\subsection{Experiment 2}
We computed the value of $a$ with $10^8$ sequences generated by the Mersenne Twister \cite{MT}.
In general, randomness tests do not use $10^8$ sequences.
However, we need to compute the value of $a$ with high accuracy to discuss the validity of (\ref{rironbunsan}) and (\ref{riron}), so we used many sequences.
Each sequence was $10^6$-bit.
The experiments were repeated 10 times and the results are presented in Table \ref{table}. 
Based on this, we can conclude that 3.7903, which was derived in the previous section, is close to the experimental result and a more accurate value compared to those generated in previous studies.

\begin{table}[h]
\caption{Value of $a$ computed using the Mersenne Twister}
\label{table}
\begin{center}
\begin{tabular}{cc}
\hline
Trial No.&$a$\\
\hline
1& 3.790208\\
2& 3.790127\\
3& 3.790568\\
4& 3.790496\\
5& 3.789058\\
6& 3.791035\\
7& 3.789565\\
8& 3.790191\\
9& 3.790522\\
10& 3.790403\\
\hline
Total &3.790217$\pm$ 0.000527\\
\hline
\end{tabular}
\end{center}
\end{table}

\section{Conclusion}
We have theoretically derived the variance of DFTT test statistics based on Parseval's theorem.
Because the variances reported by previous studies are based on numerical experiments, we can state that our result is superior subject to the tenability of the assumptions used.
In addition, and importantly, the derived value is consistent with experimental results generated herein.
Thus, we advocate that the value derived in this paper should be used when DFTT is performed so that randomness can be tested for more precisely.

\setcounter{section}{0}
\renewcommand{\appendix}{}
\appendix
\renewcommand{\thesection}{Appendix.\ \Alph{section}}
\renewcommand{\thesubsection}{\Alph{section}.\arabic{subsection}.}

\section{Proof of Proposition \ref{prop-3}}
\renewcommand{\thesection}{\Alph{section}}
We prove Proposition \ref{prop-3}.
Assume that
\begin{align*}
s_j:=\sqrt{\frac{2}{n}}\sum_{k=0}^{n-1}x_k\sin\left(\frac{2\pi kj}{n}\right),\ \ 
c_j:=\sqrt{\frac{2}{n}}\sum_{k=0}^{n-1}x_k\cos\left(\frac{2\pi kj}{n}\right).
\end{align*}
Proposition \ref{prop-3} is equivalent to the following proposition.
\begin{prop}
\label{prop-A1}
Let $R$ be an arbitrary positive integer.
As $n\to\infty$, $s_1,c_1,s_2,c_2,\cdots,s_R,c_R$ independently follow a standard normal distribution.
\end{prop}
Thus, we prove Proposition \ref{prop-A1}.
First, we introduce the following lemma.
\begin{lem}
\label{lem-A1}
Assume that $\epsilon(x)$ satisfies $\log \cos x=-\frac{1}{2}x^2+\epsilon(x)$.
Then, constants $C$ and $\Bar{x}$ exist such that
\begin{align}
|x|<\Bar{x}\ \Rightarrow\ \epsilon(x)<Cx^4.
\end{align}
\end{lem}

We define the characteristic function of the distribution followed by $(s_1,c_1,s_2,c_2,\cdots,s_R,c_R)$ as follows:
\begin{align}
\phi(u_1,v_1,u_2,v_2,\cdots,u_R,v_R):=\mathbb{E}\left[\exp\left(i\sum_{r=1}^R\left(u_rs_r+v_rc_r\right)\right)\right].
\end{align}
Then,
\begin{align*}
&\phi(u_1,v_1,u_2,v_2,\cdots,u_R,v_R)\\
=&\mathbb{E}\left[\exp\left(i\sqrt{\frac{2}{n}}\sum_{r=1}^R\left(u_r\sum_{k=0}^{n-1}x_k\sin\left(\frac{2\pi kr}{n}\right)+v_r\sum_{k=0}^{n-1}x_k\cos\left(\frac{2\pi kr}{n}\right)\right)\right)\right]\\
=&\mathbb{E}\left[\prod_{k=0}^{n-1}\exp\left(ix_k\sqrt{\frac{2}{n}}\sum_{r=1}^R\left(u_r\sin\left(\frac{2\pi kr}{n}\right)+v_r\cos\left(\frac{2\pi kr}{n}\right)\right)\right)\right]\\
=&\prod_{k=0}^{n-1}\cos\left(\sqrt{\frac{2}{n}}\sum_{r=1}^R\left(u_r\sin\left(\frac{2\pi kr}{n}\right)+v_r\cos\left(\frac{2\pi kr}{n}\right)\right)\right).
\end{align*}
For arbitrary $(u_1,v_1,u_2,v_2,\cdots,u_R,v_R)$, we have
\begin{align}
\exists D\ \text{s.t.}\ \forall k,\ \left|\sqrt{\frac{2}{n}}\sum_{r=1}^R\left(u_r\sin\left(\frac{2\pi kr}{n}\right)+v_r\cos\left(\frac{2\pi kr}{n}\right)\right)\right|\leq \frac{D}{\sqrt{n}}.
\label{1}
\end{align}
Then, if $n$ is sufficiently large,
\begin{align}
\forall k,\ \cos\left(\sqrt{\frac{2}{n}}\sum_{r=1}^R\left(u_r\sin\left(\frac{2\pi kr}{n}\right)+v_r\cos\left(\frac{2\pi kr}{n}\right)\right)\right)>0
\end{align}
for arbitrary $(u_1,v_1,u_2,v_2,\cdots,u_R,v_R)$, and so
\begin{align*}
&\log\phi(u_1,v_1,u_2,v_2,\cdots,u_R,v_R)\\
=&\sum_{k=0}^{n-1}\log\cos\left(\sqrt{\frac{2}{n}}\sum_{r=1}^R\left(u_r\sin\left(\frac{2\pi kr}{n}\right)+v_r\cos\left(\frac{2\pi kr}{n}\right)\right)\right).
\end{align*}
By Lemma \ref{lem-A1} and (\ref{1}), we have
\begin{align*}
&\log\phi(u_1,v_1,u_2,v_2,\cdots,u_R,v_R)\\
=&-\frac{1}{n}\sum_{r=1}^R\sum_{r^\prime=1}^R\left\{u_ru_{r^\prime}\sum_{k=0}^{n-1}\sin\left(\frac{2\pi kr}{n}\right)\sin\left(\frac{2\pi kr^\prime}{n}\right)\right\}\\
&\ \ \ \ -\frac{1}{n}\sum_{r=1}^R\sum_{r^\prime=1}^R\left\{v_rv_{r^\prime}\sum_{k=0}^{n-1}\cos\left(\frac{2\pi kr}{n}\right)\cos\left(\frac{2\pi kr^\prime}{n}\right)\right\}\\
&\ \ \ \ \ \ \ \ -\frac{2}{n}\sum_{r=1}^R\sum_{r^\prime=1}^R\left\{u_rv_{r^\prime}\sum_{k=0}^{n-1}\sin\left(\frac{2\pi kr}{n}\right)\cos\left(\frac{2\pi kr^\prime}{n}\right)\right\}\\
&\ \ \ \ \ \ \ \ \ \ \ \ +\sum_{k=0}^{n-1}\epsilon\left(\sqrt{\frac{2}{n}}\sum_{r=1}^R\left(u_r\sin\left(\frac{2\pi kr}{n}\right)+v_r\cos\left(\frac{2\pi kr}{n}\right)\right)\right)\\
=&-\frac{1}{n}\sum_{r=1}^R\left\{u_r^2\sum_{k=0}^{n-1}\sin^2\left(\frac{2\pi kr}{n}\right)\right\}-\frac{1}{n}\sum_{r=1}^R\left\{v_r^2\sum_{k=0}^{n-1}\cos^2\left(\frac{2\pi kr}{n}\right)\right\}+O\left(\frac{1}{n}\right)\\
=&\ \ \sum_{r=1}^R\left(-\frac{1}{2}u_r^2-\frac{1}{2}v_r^2\right)+O\left(\frac{1}{n}\right).
\end{align*}
For arbitrary $(u_1,v_1,u_2,v_2,\cdots,u_R,v_R)$,
\begin{align*}
\lim_{n\to\infty}\phi(u_1,v_1,u_2,v_2,\cdots,u_R,v_R)=&\exp\left(\sum_{r=1}^R\left(-\frac{1}{2}u_r^2-\frac{1}{2}v_r^2\right)\right)\\
=&\prod_{r=1}^R\exp\left(-\frac{1}{2}u_r^2\right)\exp\left(-\frac{1}{2}v_r^2\right).
\end{align*}
Thus, we have proved Proposition \ref{prop-A1}.
 Then, it follows from Proposition \ref{prop-A1} that Proposition \ref{prop-3} holds.

\end{document}